\documentclass[12pt,twoside]{article}

\usepackage{amsmath,amssymb,amsthm}

\newcommand{\Prim}{\mathrm{Prim}}
\newcommand{\vol}{\mathrm{vol}}
\newcommand{\Hyp}{\mathrm{Hyp}}

\newcommand{\Tr}{\mathrm{Tr}}
\newcommand{\tr}{\mathrm{tr}}
\newcommand{\Ind}{\mathrm{Ind}}

\newcommand{\as}{\quad\text{as}\quad}
\newcommand{\tinf}{\to\infty}
\newcommand{\disp}{\displaystyle}
\newcommand{\bsla}{\backslash}
\newcommand{\nt}{\notag}

\newcommand{\D}{\mathfrak{D}}
\newcommand{\M}{\mathcal{M}}
\newcommand{\bR}{\mathbb{R}}
\newcommand{\bZ}{\mathbb{Z}}
\newcommand{\sr}{\mathrm{SL}_2(\bR)}
\newcommand{\sz}{\mathrm{SL}_2(\bZ)}
\newcommand{\divset}{\hspace{3pt}|\hspace{3pt}}
\newcommand{\bigdivset}{\hspace{3pt}\big|\hspace{3pt}}
\newcommand{\Bigdivset}{\hspace{3pt}\Big|\hspace{3pt}}

\newcommand{\mmid}{\hspace{2pt}||\hspace{2pt}}

\newtheorem{thm}{Theorem}[section]

\newtheorem{lem}[thm]{Lemma}
\newtheorem{cor}[thm]{Corollary}
\newtheorem{rem}[thm]{Remark}

\numberwithin{equation}{section}

\title{Arithmetic expressions of Selberg's zeta functions for congruence subgroups}
\author{Yasufumi Hashimoto}
\date{}
\begin{document}
\maketitle 
\begin{abstract}
In \cite{Sa}, it was proved that the Selberg zeta function for $\sz$ 
is expressed in terms of the fundamental units and the class numbers 
of the primitive indefinite binary quadratic forms.
The aim of this paper is to obtain similar arithmetic expressions 
of the logarithmic derivatives of the Selberg zeta functions for congruence subgroups of $\sz$.
As applications, we study the Brun-Titchmarsh type prime geodesic theorem, 
the asymptotic behavior of the sum of the class number.
\end{abstract}

\renewcommand{\thefootnote}{}

\section{Introduction}


Let $H$ be the upper half plane and 
$\Gamma$ a discrete subgroup of $\sr$ such that $\vol(\Gamma\bsla H)<\infty$. 
We define the Selberg zeta function by
\begin{align}
Z_{\Gamma}(s):=\prod_{p\in\Prim(\Gamma)}\prod_{n=0}^{\infty}\big(1-N(p)^{-s-n}\big)\qquad\Re{s}>1,
\end{align}
where $\Prim(\Gamma)$ is the set of the primitive hyperbolic conjugacy classes of $\Gamma$ 
and $N(p)$ is a square of the larger eigenvalue of $p\in\Prim(\Gamma)$.
In \cite{Sa}, it was shown that $Z_{\Gamma}(s)$ and $Z'_{\Gamma}(s)/Z_{\Gamma}(s)$ for $\Gamma=\sz$ have
the following expression.
\begin{align}
Z_{\Gamma}(s)=&\prod_{D\in\D}\prod_{n=0}^{\infty}\big(1-\epsilon(D)^{-2(s+n)}\big)^{h(D)},\label{prodselb}\\
\frac{Z'_{\Gamma}(s)}{Z_{\Gamma}(s)}
=&\sum_{D\in\D}
\sum_{j=1}^{\infty}h(D)\frac{2\log{\epsilon(D)}}{1-\epsilon(D)^{-2j}}\epsilon(D)^{-2js}\label{logderiselb},
\end{align}
where $\D:=\{D>0|D\equiv0,1\bmod 4,\text{not a square}\}$, $\epsilon(D)$ is the fundamental unit and  
$h(D)$ is the class number of discriminant $D\in\D$ in the narrow sense.
In \cite{AKN}, analogues of the expressions \eqref{prodselb} and \eqref{logderiselb} are also established in the case 
where $\Gamma$ is a quaternion group. 

In this paper, we study the arithmetic expressions of the Selberg zeta functions 
for congruence subgroups of $\sz$ defined by
\begin{align*}
\Gamma_0(N)=&\Big\{\gamma=\big(\gamma_{ij}\big)_{1\leq i,j\leq2}\in \sz\Bigdivset
\gamma_{21}\equiv0\bmod N\Big\},\\
\Gamma_1(N)=&\Big\{\gamma=\big(\gamma_{ij}\big)_{1\leq i,j\leq2}\in \sz\Bigdivset
\gamma_{11}\equiv\gamma_{22}\equiv\pm1,\gamma_{21}\equiv0\bmod N\Big\},\\
\Gamma(N)=&\Big\{\gamma=\big(\gamma_{ij}\big)_{1\leq i,j\leq2}\in \sz\Bigdivset
\gamma_{11}\equiv\gamma_{22}\equiv\pm1,\gamma_{21}\equiv\gamma_{12}\equiv0\bmod N\Big\}.
\end{align*}
Note that 
\begin{align*}
[\sz:\Gamma_0(N)]=&N\prod_{p|N}(1+p^{-1}),\\
[\sz:\Gamma_1(N)]=&\begin{cases}3&(N=2),\\ \disp\frac{N^2}{2}\prod_{p|N}(1-p^{-2}),&(N\geq3),
\end{cases}\\
[\sz:\Gamma(N)]=&\begin{cases}6&(N=2),\\ \disp\frac{N^3}{2}\prod_{p|N}(1-p^{-2})&(N\geq3).
\end{cases}
\end{align*}
The main results of the present paper is as follows. 
\begin{thm}\label{thm1}
Let $N>1$ be an integer and $N=\prod_{p\mid N}p^r$ the factorization of $N$. 
For $D\in\D$, denote by $(t_j(D),u_j(D))$ the $j$-th positive solution of the Pell equation $t^2-u^2D=4$.
Then the logarithmic derivative of the Selberg zeta functions 
for $\Gamma=\Gamma_0(N)$, $\Gamma_1(N)$ or $\Gamma(N)$ are expressed as
\begin{align}\label{congselb}
\frac{Z'_{\Gamma}(s)}{Z_{\Gamma}(s)}
=&\sum_{D\in\D}
\sum_{j=1}^{\infty}M_{\Gamma}(D,j)
h(D)\frac{2\log{\epsilon(D)}}{1-\epsilon(D)^{-2j}}\epsilon(D)^{-2js},
\end{align}
where $M_{\Gamma}(D,j)$'s are given by
\begin{align*}
M_{\Gamma_0(N)}(D,j)=&\prod_{p\mid N}L^{(0)}_{p^r}(D,j),\\
M_{\Gamma_1(N)}(D,j)=&\begin{cases}\disp\frac{1}{2}\prod_{p\mid N}L^{(1)}_{p^r}(D,j)&(t_j(D)\equiv\pm2\bmod N),\\
0&(\text{otherwise}),\end{cases}\\
M_{\Gamma(N)}(D,j)=&\begin{cases}\disp\frac{N^3}{2}\prod_{p\mid N}(1-p^{-2})
&(N>2,t_j(D)\equiv\pm2\bmod N,N\mid u_j(D)),\\
6&(N=2,2\mid u_j(D)),\\
0&(\text{otherwise}).\end{cases}
\end{align*}
Here $L^{(0)}_{p^r}(D,j)$ is defined as follows; when $p^r=2$
\begin{align*}
L^{(0)}_{2}(D,j):=&\begin{cases}3&(2\mid u_j(D)),\\ 
1&(2\nmid u_j(D), 4\mid D),\\
0&(\text{otherwise}).\end{cases}
\intertext{when $p^r=2^r\geq4$;}
L^{(0)}_{2^r}(D,j):=&\begin{cases}3\times 2^{r-1}&(2^r\mid u_j(D)),\\ 
2^{r-1}&(\text{$2^{r-1}\mmid u_j(D)$ or $2^{r-2}\mmid u_j(D),8\mid D$}),\\ 
2^{[\frac{r+k}{2}]}&(2^{k}\mmid u_j(D),2^{r-k}\mid D,k\leq r-3),\\
2^{m+k+2}&\big(2^{k}\mmid u_j(D),2^{2m}\mid D,D/2^{2m}\equiv1\bmod{8},\\&1< 2m<r-k-2),\\
2^{m+k+1}&\big(2^{k}\mmid u_j(D),2^{2m}\mid D,D/2^{2m}\equiv1\bmod{8},\\&m=[\frac{r-k-1}{2}]),\\
0&(\text{otherwise}),
\end{cases}\\
\intertext{when $p\geq 3$;}
L^{(0)}_{p^r}(D,j):=&\begin{cases}p^{r-1}(p+1)&(p^r\mid u_j(D)),\\ 
p^{[\frac{r+k}{2}]}&(p^k\mmid u_j(D),p^{r-k}|D,k\leq r-1),\\
2p^{m+k}&\big(p^k\mmid u_j(D),p^{2m}|D,\big(\frac{D/p^{2m}}{p}\big)=1,
\\&0\leq 2m<r-k),\\
0&(\text{otherwise}),
\end{cases}
\end{align*}
where $p^k\mmid a$ means that $p^k\mid a$ and $p^{k+1}\nmid a$.
On the other hand, $L^{(1)}_{p^r}(D,j)$ is given by
\begin{align*}
L^{(1)}_{p^r}(D,j):=&\begin{cases}p^{2r-2}(p^2-1)&(p^r\mid u_j(D)),\\ 
p^{r+k-1}(p-1)&(p^k\mmid u_j(D),p^{r-k}|D,k\leq r-1).\end{cases}
\end{align*}
\end{thm}

As a corollary of Theorem \ref{thm1}, we estimate the number of primitive hyperbolic conjugacy classes $p$ of $\Gamma$
such that $x\leq N(p)<x+y$.
Recall the prime geodesic theorem;
\begin{align}\label{prim}
\pi_{\Gamma}(x)\sim \frac{x}{\log{x}}\as x\tinf,
\end{align}
where $\pi_{\Gamma}(x)$ is the number of primitive hyperbolic conjugacy classes of $\Gamma$ 
whose norm is smaller than $x$.
In \cite{Iw} and \cite{AKN}, when $\Gamma$ is $\sz$ or a quaternion group, the following estimate is shown   
for $x^{1/2}(\log{x})^2<y<x$.
\begin{align}\label{brun}
\pi_{\Gamma}(x+y)-\pi_{\Gamma}(x)\ll y.
\end{align} 
In this paper, applying Theorem \ref{thm1} with the help of the result in \cite{Iw}, 
we also prove the following estimate 
in the case where $\Gamma=\Gamma_0(N)$, $\Gamma_1(N)$ or $\Gamma(N)$. 
\begin{cor}\label{thm2}
If $\Gamma=\Gamma_0(N)$, $\Gamma_1(N)$ or $\Gamma(N)$, then, for $x^{1/2}(\log{x})^2<y<x$, we have
\begin{align}
\pi_{\Gamma}(x+y)-\pi_{\Gamma}(x)\ll y,
\end{align}
where the implied constant depends only on $\Gamma$.
\end{cor}

Theorem \ref{thm1} has another corollary related with the prime geodesic theorem.
Owing to \eqref{logderiselb}, 
the following estimate yields from the prime geodesic theorem for $\Gamma=\sz$.
\begin{align}\label{class}
\sum_{\begin{subarray}{c}D\in\D\\\epsilon(D)<x\end{subarray}}&h(D)\sim \frac{x^2}{2\log{x}}\as x\tinf.
\end{align}
As a development of \eqref{class}, Sarnak \cite{Sa} have estimated the asymptotic behavior of 
the sum of the class number $h(D)$ over $D\in\D$ such that $p\mid u_1$ and $\epsilon(D)<x$. 
In this paper, by using the result of Theorem \ref{thm1} for $\Gamma_1(p)$ and $\Gamma(p)$, 
we give another development of \eqref{class} as follows.
\begin{cor}\label{cor3}
Let $p\geq3$ be a prime number. Then we have 
\begin{align}
\sum_{\begin{subarray}{c}D\in\D,p\mid D \\ \epsilon(D)<x\end{subarray}}&h(D)\sim C(p)\frac{x^2}{\log{x}}\as x\tinf,
\end{align}
where $C(p)$ is a constant such that $1/p\leq C(p)\leq p/(p^2-1)$.
\end{cor}

\section{Proofs of Theorem \ref{thm1} and corollaries}
\subsection{Lemmas}
Before proving the theorem, we prepare the following lemmas.
First, we recall the following Venkov-Zograf's theorem (see \cite{VZ}).
\begin{lem}\label{venkov}
Let $\Gamma$ be a discrete subgroup of $\sr$ such that 
$\vol(\Gamma\bsla H)<\infty$, 
$\Gamma'$ a subgroup of $\Gamma$ with a finite index in $\Gamma$, 
$\chi$ a finite dimensional unitary representation of $\Gamma'$, and 
\begin{align*}
Z_{\Gamma}(s,\chi):=\prod_{p\in\Prim(\Gamma)}\prod_{n=0}^{\infty}\det\big(I-\chi(p)N(p)^{-s-n}\big).
\end{align*}
Then, we have 
\begin{align}\label{induced}
Z_{\Gamma'}(s,\chi)=Z_{\Gamma}\big(s,\Ind_{\Gamma'}^{\Gamma}\chi\big).
\end{align} 
\end{lem}

\begin{lem}\label{lem1}
(1) Let $p$ be a prime number.
Then the complete system of representatives of $\sz/\Gamma_0(p^r)$ can be chosen by 
\begin{align*}
\bigg\{\begin{pmatrix}1&0\\m&1\end{pmatrix},\quad \begin{pmatrix}lp&-1\\1&0\end{pmatrix}\Bigdivset
 m\in\bZ/p^r\bZ,l\in\bZ/p^{r-1}\bZ\bigg\}.
\end{align*}
(2) Let $p$ be a prime number and $N>1$ an integer which is relatively prime to $p$. 
Then the complete system of representatives of $\Gamma_0(N)/\Gamma_0(Np^r)$ can be chosen by 
\begin{align*}
\bigg\{\begin{pmatrix}1&0\\mN&1\end{pmatrix},\quad \begin{pmatrix}lp& k_1\\N&k_2\end{pmatrix}
\Bigdivset m\in\bZ/p^r\bZ,l\in\bZ/p^{r-1}\bZ\bigg\},
\end{align*}
where $k_1,k_2$ are integers such that $lpk_2-k_1N=1$.\\
\end{lem}
The proof of the lemma above is elementary. Then we omit the proof.
Now, for convenience, we put $t_{\gamma}:=\tr{\gamma}$, 
$u_{\gamma}:=\gcd(\gamma_{21},\gamma_{12},\gamma_{22}-\gamma_{11})_{>0}$ and 
$D_{\gamma}:=(t_{\gamma}^2-4)/u_{\gamma}^2$.
We can prove the following lemma by using Lemma \ref{lem1}.
\begin{lem}\label{lem2}
(1) Let $p$ be a prime number and $\gamma\in \sz$. Then we have 
\begin{align}
\Tr\Big(\big(\Ind_{\Gamma_0(p^r)}^{\sz}1\big)(\gamma)\Big)=\hat{L}^{(0)}_{p^r}(t_{\gamma},u_{\gamma}),
\end{align}
where $d_{t,u}:=(t^2-4)/u^2$ and $\hat{L}^{(0)}_{p^r}(t,u)$ is defined as follows; when $p^r=2$
\begin{align*}
\hat{L}^{(0)}_{2}(t,u):=&\begin{cases}3&(2\mid u),\\ 
1&(2\nmid u, 4\mid d_{t,u}),\\
0&(\text{otherwise}).\end{cases}
\intertext{when $p^r=2^r\geq4$;}
\hat{L}^{(0)}_{2^r}(t,u):=&\begin{cases}3\times 2^{r-1}&(2^r\mid u),\\ 
2^{r-1}&(\text{$2^{r-1}\mmid u$ or $2^{r-2}\mmid u,8\mid d_{t,u}$}),\\ 
2^{[\frac{r+k}{2}]}&(2^{k}\mmid u,2^{r-k}\mid d_{t,u},k\leq r-3),\\
2^{m+k+2}&\big(2^{k}\mmid u,2^{2m}\mid d_{t,u},d_{t,u}/2^{2m}\equiv1\bmod{8},\\&1<2m<r-k-2),\\
2^{m+k+1}&\big(2^{k}\mmid u,2^{2m}\mid d_{t,u},d_{t,u}/2^{2m}\equiv1\bmod{8},\\&m=[\frac{r-k-1}{2}]),\\
0&(\text{otherwise}),
\end{cases}\\
\intertext{when $p\geq 3$;}
\hat{L}^{(0)}_{p^r}(t,u):=&\begin{cases}p^{r-1}(p+1)&(p^r\mid u),\\ 
p^{[\frac{r+k}{2}]}&(p^k\mmid u,p^{r-k}\mid d_{t,u},k\leq r-1),\\
2p^{m+k}&\big(p^k\mmid u,p^{2m}\mid d_{t,u},\big(\frac{d_{t,u}/p^{2m}}{p}\big)=1,
\\&0\leq 2m<r-k\big),\\
0&(\text{otherwise}).
\end{cases}
\end{align*}

\noindent(2) Let $p$ be a prime number, $N>1$ an integer 
which is relatively prime to $p$ and $\gamma\in\Gamma_0(N)$. 
Then we have
\begin{align}
\Tr\Big(\big(\Ind_{\Gamma_0(Np^r)}^{\Gamma_0(N)}1\big)(\gamma)\Big)=\hat{L}^{(0)}_{p^r}(t_{\gamma},u_{\gamma}).
\end{align}
\end{lem}

\begin{proof}
Let $\{a_k\}_{1\leq k\leq p^{r-1}(p+1)}$ be the complete system of representatives of $\sz/\Gamma_0(p^r)$. 
The induced representation $\Ind_{\Gamma_0(p^r)}^{\sz}1$ is written as
\begin{align*}
\big(\Ind_{\Gamma_0(p^r)}^{\sz}1\big)(\gamma)=\Big(\delta(a_i^{-1}\gamma a_j)\Big)_{1\leq i,j\leq p^{r-1}(p+1)},
\end{align*}
where 
\begin{align*}
\delta(g)=\begin{cases}1&\text{if}\quad g\in\Gamma_0(p^r),\\
0&\text{if}\quad g\not\in\Gamma_0(p^r).\end{cases}
\end{align*}
Hence we have
\begin{align*}
\Tr\Big(\big(\Ind_{\Gamma_0(p^r)}^{\sz}1\big)(\gamma)\Big)
=&\#\Big\{1\leq i\leq p^{r-1}(p+1)\bigdivset a_i^{-1}\gamma a_i\in\Gamma_0(p^r)\Big\}\\
=&\#\Big\{1\leq i\leq p^{r-1}(p+1)\bigdivset (a_i^{-1}\gamma a_i)_{21}\equiv0\bmod p\Big\}.
\end{align*}
According to (1) of Lemma \ref{lem1}, we have
\begin{align}\label{case}
\Tr\Big(\big(\Ind_{\Gamma_0(p^r)}^{\sz}1\big)(\gamma)\Big)
=&\#\{m\in\bZ/p^r\bZ\bigdivset \gamma_{12}m^2+(\gamma_{11}-\gamma_{22})m-\gamma_{21}\equiv0\bmod{p^r}\}\nt\\
+&\#\{l\in\bZ/p^{r-1}\bZ\bigdivset p^2\gamma_{21}l^2+p(\gamma_{11}-\gamma_{22})l-\gamma_{12}\equiv0\bmod{p^r}\}
\end{align}

\noindent {\bf (i) The case $p^r\mid u_{\gamma}$:} \\
Since $\gamma_{21},\gamma_{12},\gamma_{11}-\gamma_{22}\equiv0\bmod p^r$, 
$(a_i^{-1}\gamma a_i)_{21}\equiv0$ holds for any $i$. 
Hence we have 
\begin{align*}
\Tr\Big(\big(\Ind_{\Gamma_0(p^r)}^{\sz}1\big)(\gamma)\Big)=p^{r-1}(p+1).
\end{align*}
\noindent{\bf (ii) The case $p^k\mmid u_{\gamma}$ $(k\leq r-1)$:}\\
Let $A:=\gamma_{12}/p^k$, $B:=(\gamma_{11}-\gamma_{22})/p^k$ and $C:=-\gamma_{21}/p^k$. 
Then we have 
\begin{align*}
\Tr\Big(\big(\Ind_{\Gamma_0(p^r)}^{\sz}1\big)(\gamma)\Big)
=&\#\{m\in\bZ/p^r\bZ\divset Am^2+Bm+C\equiv0\bmod{p^{r-k}}\}\\
+&\#\{l\in\bZ/p^{r-1}\bZ\divset -p^2Cl^2+pBl-A\equiv0\bmod{p^{r-k}}\}.
\end{align*}
For simplicity, we treat only the case where $p$ is odd and divides neither $A$, $B$ nor $C$. 
In this case, it is easy to see that 
\begin{align}
\Tr\Big(\big(\Ind_{\Gamma_0(p^r)}^{\sz}1\big)(\gamma)\Big)
=&\#\{m\in(\bZ/p^r\bZ)^{*}|Am^2+Bm+C\equiv0\bmod{p^{r-k}}\}.\label{equiv1}
\end{align}
Since $A$ is not divided by $p$, we can rewrite the equation in \eqref{equiv1} as follows.
\begin{align}
(m+(2A)^{-1}B)^2\equiv (2A)^{-2}(B^2+4AC)\equiv (2A)^{-2}D_{\gamma}\bmod{p^{r-k}}.\label{equiv2}
\end{align}
Counting the number of solutions of \eqref{equiv2}, 
we see that $\Tr\Big(\big(\Ind_{\Gamma_0(p^r)}^{\sz}1\big)(\gamma)\Big)$ 
coincides $\hat{L}^{(0)}_{p^r}(t_{\gamma},u_{\gamma})$. 

In any other cases ($p=2$ or $p$ divides at least one of $A,B,C$ and does not divide at least one of them),
the results are contained in that of the case above.
Hence the proof of (1) is completed. 
We can also prove (2) similarly. 
\end{proof}

\begin{lem}\label{lem4}
Let $t\geq3$ and $u\geq1$ be integers such that $d_{t,u}\in\D$.
Then we have
\begin{align}
\#\{\gamma\in\Hyp(\sz)\divset t_{\gamma}=t,u_{\gamma}=u\}=h(d_{t,u}).
\end{align}
\end{lem}
\begin{proof} See, for example, [Sa]. \end{proof}

\subsection{Proof of Theorem \ref{thm1}}
The expression \eqref{congselb} can be rewritten as 
\begin{align}\label{tarith}
\frac{Z'_{\Gamma}(s)}{Z_{\Gamma}(s)}=\sum_{t=3}^{\infty}
\sum_{u\in U(t)}
\hat{M}_{\Gamma}(t,u)
h\big(d_{t,u}\big)j(t,u)^{-1}\frac{2\log{\epsilon(t)}}{1-\epsilon(t)^{-2}}\epsilon(t)^{-2s},
\end{align}
where 
\begin{align*}
&U(t):=\{u\geq1| d_{t,u}\in\D\},\\
&\epsilon(t)=\frac{1}{2}\bigl(t+\sqrt{t^2-4}\bigr)=\frac{1}{2}\bigl(t+u\sqrt{d_{t,u}}\bigr),\\
&j(t,u)=\max\Bigl\{j\geq1 \Big|\epsilon(t)= 
\Bigl(\frac{1}{2}\big(t_0+u_0\sqrt{d_{t,u}}\big)\Bigr)^j,\quad\exists t_0,u_0\geq 1\Bigr\},
\end{align*}
and $\hat{M}_{\Gamma}(t,u)$ is given by
\begin{align*}
\hat{M}_{\Gamma_0(N)}(t,u)=&\prod_{p\mid N}\hat{L}^{(0)}_{p^r}(t,u),\\
\hat{M}_{\Gamma_1(N)}(t,u)=&\begin{cases}\disp\frac{1}{2}\prod_{p\mid N}\hat{L}^{(1)}_{p^r}(t,u)&(t\equiv\pm2\bmod N),\\
0&(\text{otherwise}),\end{cases}\\
\hat{M}_{\Gamma(N)}(t,u)=&\begin{cases}\disp\frac{N^3}{2}\prod_{p\mid N}(1-p^{-2})&(t\equiv\pm2\bmod N,N\mid u),\\
0&\text{(otherwise)}.\end{cases}
\end{align*}
Here, 
$\hat{L}^{(0)}_{p^r}(t,u)$ is defined in Lemma \ref{lem2} and
\begin{align*}
&\hat{L}^{(1)}_p(t,u):=\begin{cases}p^{2r-2}(p^2-1)&(p^r\mid u),\\ 
p^{r+k-1}(p-1)&(p^k\mmid u,p^{r-k}\mid d_{t,u},k\leq r-1).\end{cases}
\end{align*}
Hence, we prove the expression \eqref{tarith} instead of \eqref{congselb}.

Let $N>1$ be an integer, and $\Gamma=\Gamma_0(N)$, $\Gamma_1(N)$ or $\Gamma(N)$.
According to Lemma \ref{venkov}, we have 
\begin{align}\label{selb}
Z_{\Gamma}(s)=Z_{\sz}\Big(s,\Ind_{\Gamma}^{\sz}1\Big).
\end{align}
The logarithmic derivative of \eqref{selb} is written by
\begin{align}
\frac{Z'_{\Gamma}(s)}{Z_{\Gamma}(s)}=\sum_{\gamma\in\Hyp(\sz)}
\Tr{\Big(\big(\Ind_{\Gamma}^{\sz}1\big)(\gamma)\Big)}
\frac{j_{\gamma}^{-1}\log{N(\gamma)}}{1-N(\gamma)^{-1}}N(\gamma)^{-s},\label{logderi}
\end{align}
where $j_{\gamma}\geq1$ is an integer such that $\gamma=\delta^{j_{\gamma}}$ 
for some $\delta\in\Prim(\Gamma)$.
The right hand side of \eqref{logderi} is expressed as follows.
\begin{align}
\frac{Z'_{\Gamma}(s)}{Z_{\Gamma}(s)}
=&\sum_{t=3}^{\infty}\frac{2\log{\epsilon(t)}}{1-\epsilon(t)^{-2}}\epsilon(t)^{-2s}
\sum_{\gamma\in A(t)}j_{\gamma}^{-1}\Tr{\Big(\big(\Ind_{\Gamma}^{\sz}1\big)(\gamma)\Big)}\notag\\
=&\sum_{t=3}^{\infty}\frac{2\log{\epsilon(t)}}{1-\epsilon(t)^{-2}}\epsilon(t)^{-2s}
\sum_{u\in U(t)}j(t,u)^{-1}
\sum_{\gamma\in B(t,u)}\Tr{\Big(\big(\Ind_{\Gamma}^{\sz}1\big)(\gamma)\Big)},\label{logderi2}
\end{align}
where 
\begin{align*}
A(t):=&\{\gamma\in\Hyp(\sz)\divset t_{\gamma}=t\},\\
B(t,u):=&\{\gamma\in\Hyp(\sz)\divset t_{\gamma}=t,u_{\gamma}=u\}.
\end{align*}

Now, we calculate $\Tr{\Big(\big(\Ind_{\Gamma}^{\sz}1\big)(\gamma)\Big)}$ for each $\Gamma$. 

\noindent{\bf 1) The case $\Gamma=\Gamma_0(N)$:}\\
Let $p$ be a prime number which divides $N$ 
and $\{a_k\}$ (resp. $\{b_l\}$) a complete set of representatives of $\Gamma_0(N/p^r)/\Gamma_0(N)$ 
(resp. $\sz/\Gamma_0(N/p^r)$).
The induced representation $\big(\Ind_{\Gamma_0(N)}^{\sz}1\big)(\gamma)$ is written as
\begin{align}
\big(\Ind_{\Gamma_0(N)}^{\sz}1\big)(\gamma)
=&\Big(\Ind_{\Gamma_0(N/p^r)}^{\sz}\big(\Ind_{\Gamma_0(N)}^{\Gamma_0(N/p^r)}1\big)\Big)(\gamma)\notag\\
=&\Big(\delta'(b_i^{-1}\gamma b_j)\Big)_{1\leq i,j\leq [\sz:\Gamma_0(N/p^r)]},\label{ind}
\end{align}
where 
\begin{align*}
\delta'(g)=\begin{cases}\big(\Ind_{\Gamma_0(N)}^{\Gamma_0(N/p^r)}1\big)(g)&\text{if $g\in\Gamma_0(N/p^r)$},\\
0&\text{if $g\notin\Gamma_0(N/p^r)$}.\end{cases}
\end{align*}
Hence, the trace of \eqref{ind} is as follows. 
\begin{align*}
\Tr\Big(\big(\Ind_{\Gamma_0(N)}^{\sz}1\big)(\gamma)\Big)
=\sum_{i\in\M}\Tr\Big(\big(\Ind_{\Gamma_0(N)}^{\Gamma_0(N/p^r)}1\big)(b_i^{-1}\gamma b_i)\Big),
\end{align*}
where
\begin{align*}
\M:=\Big\{1\leq i\leq [\sz:\Gamma_0(N/p^r)]\bigdivset b_i^{-1}\gamma b_i\in\Gamma_0(N/p^r)\Big\}.
\end{align*}
By using (2) of Lemma \ref{lem2} and the fact that 
$(t_{g^{-1}\gamma g},u_{g^{-1}\gamma g})=(t_{\gamma},u_{\gamma})$ for any $g\in \sz$, we have
\begin{align*}
\Tr\Big(\big(\Ind_{\Gamma_0(N)}^{\sz}1\big)(\gamma)\Big)
=&\sum_{i\in\M}\hat{L}^{(0)}_p(t_{b_i^{-1}\gamma b_i},u_{b_i^{-1}\gamma b_i})\\
=&\sum_{i\in\M}L^{(0)}_p(t_{\gamma},u_{\gamma})\\
=&\hat{L}^{(0)}_p(t_{\gamma},u_{\gamma})\Tr\Big(\big(\Ind_{\Gamma_0(N/p^r)}^{\sz}1\big)(\gamma)\Big).
\end{align*}
This means that
\begin{align*}
\Tr\Big(\big(\Ind_{\Gamma_0(N)}^{\sz}1\big)(\gamma)\Big)
=\prod_{p\mid N}\hat{L}^{(0)}_p(t_{\gamma},u_{\gamma})
\end{align*}
holds recursively.
Hence, owing to Lemma \ref{lem4}, we have
\begin{align}\label{number}
\sum_{\gamma\in B(t,u)}\Tr\Big(\big(\Ind_{\Gamma_0(N)}^{\sz}1\big)(\gamma)\Big)
=\bigg(\prod_{p\mid N}\hat{L}^{(0)}_p(t,u)\bigg)\sum_{\gamma\in B(t,u)}1=\hat{M}_{\Gamma_0(N)}(t,u)h(d_{t,u}).
\end{align}
Therefore, substituting \eqref{number} to \eqref{logderi2}, 
we obtain the expression \eqref{tarith} in the case of $\Gamma_0(N)$.

\noindent{\bf 2) The case $\Gamma=\Gamma_1(N)$:}\\
It is easy to see that $\Gamma_1(N)$ is a normal subgroup of $\Gamma_0(N)$. 
Then we have
\begin{align*}
\Tr\Big(\big(\Ind_{\Gamma_1(N)}^{\sz}1\big)(\gamma)\Big)
=&\Tr\Big(\Big(\Ind_{\Gamma_0(N)}^{\sz}\big(\Ind_{\Gamma_1(N)}^{\Gamma_0(N)}1\big)\Big)(\gamma)\Big)\\
=&\sum_{\begin{subarray}{c}g\in \sz/\Gamma_0(N)\\ g^{-1}\gamma g\in\Gamma_0(N)\end{subarray}}
\Tr\Big(\big(\Ind_{\Gamma_1(N)}^{\Gamma_0(N)}1\big)(g^{-1}\gamma g)\Big)\\
=&\sum_{\begin{subarray}{c}g\in \sz/\Gamma_0(N)\\ g^{-1}\gamma g\in\Gamma_1(N)\end{subarray}}[\Gamma_0(N):\Gamma_1(N)]\\
=&[\Gamma_0(N):\Gamma_1(N)]\times\#\big\{g\in \sz/\Gamma_0(N)\divset g^{-1}\gamma g\in\Gamma_1(N)\big\}.
\end{align*}
Hence, by using (1) and (2) of Lemma \ref{lem1}, we can obtain the desired result similar to the case of $\Gamma_0(N)$.

\noindent{\bf 3) The case $\Gamma=\Gamma(N)$:}\\
Since $\Gamma(N)$ is a normal subgroup of $\sz$, it is easy to see that 
\begin{align*}
\big(\Ind_{\Gamma_1(N)}^{\sz}1\big)(\gamma)
=&\begin{cases}[\sz:\Gamma(N)]&(\gamma\in \Gamma(N)),\\0&(\text{otherwise}).\end{cases}
\end{align*}
The condition $\gamma\in\Gamma(N)$ is equivalent to $t_j(D)\equiv\pm2\pmod{N}$ and $N\mid u_j(D)$.
Hence \eqref{tarith} also holds for $\Gamma=\Gamma(N)$.

This completes the proof. 
\qed

\subsection{Proof of Corollary \ref{thm2}}
Put
\begin{align*}
\hat{\pi}_{\Gamma}(x):=\sum_{\begin{subarray}{c} \gamma\in\Prim(\Gamma),j\geq1\\ N(\gamma)^j<x\end{subarray}}j^{-1}.
\end{align*} 
Because of Theorem \ref{thm1}, $\hat{\pi}_{\Gamma}(x)$ is expressed as
\begin{align*}
\hat{\pi}_{\Gamma}(x)
=&\sum_{\begin{subarray}{c} D\in\D,j\geq1\\ \epsilon(D)^{2j}<x\end{subarray}}j^{-1}M_{\Gamma}(D,j)h(D).
\end{align*}
Since $0\leq M_{\Gamma}(D,j)\leq [\sz:\Gamma]$ for any $D\in\D, j\geq1$, we have
\begin{align}
\hat{\pi}_{\Gamma}(x+y)-\hat{\pi}_{\Gamma}(x)
\leq&[\sz:\Gamma]
\sum_{\begin{subarray}{c} D\in\D,j\geq1\\ x\leq\epsilon(D)^{2j}<x+y\end{subarray}}j^{-1}h(D)\notag\\
=&[\sz:\Gamma]\big(\hat{\pi}_{\sz}(x+y)-\hat{\pi}_{\sz}(x)\big).\label{sl2}
\end{align}
The following estimate was shown in \cite{Iw}.
\begin{align}\label{iwaniec}
\pi_{\sz}(x+y)-\pi_{\sz}(x)\ll y.
\end{align}
Furthermore, on account of $\hat{\pi}_{\Gamma}(x)=\sum_{j\geq1}j^{-1}\pi_{\Gamma}(x^{1/j})$,
it is easy to see that 
\begin{align}
&0\leq\hat{\pi}_{\Gamma}(x)-\pi_{\Gamma}(x)\ll x^{1/2}\ll y.\label{prim1}
\end{align}
Hence, combining \eqref{sl2}, \eqref{iwaniec} and \eqref{prim1}, we can obtain the desired result. 
\qed

\subsection{Proof of Corollary \ref{cor3}}
Let $p$ be an odd prime number and $j$ a positive integer. 
We define the following subsets of $\D$. 
\begin{align*}
\D_{p}(x):=&\{D\in\D \divset p\mid D,\epsilon(D)<x\},\\
\D_{p,j}^{(1)}(x):=&\{D\in\D \divset \text{$p\mid D$ or $p\mid u_j(D)$},\epsilon(D)^j<x\},\\
\D_{p,j}^{(2)}(x):=&\{D\in\D \divset p\mid D,p\nmid u_j(D),\epsilon(D)^j<x\}.
\end{align*}
Since $\D_{p,1}^{(2)}(x)\subset\D_{p}(x)\subset\D_{p,1}^{(1)}(x)$, 
we should get the following asymptotic behaviors to prove Corollary \ref{cor3}. 
\begin{align}
\sum_{D\in\D_{p,1}^{(1)}(x)}h(D)\sim&\frac{p}{p^2-1}\frac{x^2}{\log{x}},
\label{d1}\\
\sum_{D\in\D_{p,1}^{(2)}(x)}h(D)\sim&\frac{1}{p}\frac{x^2}{\log{x}}
\as x\tinf.\label{d2}
\end{align}

According to Theorem \ref{thm1}, we can write as
\begin{align*}
\hat{\pi}_{\Gamma_1(p)}(x^2)
=&\frac{1}{2}(p^2-1)\sum_{j\geq1}
\sum_{\begin{subarray}{c}D\in\D,p\mid u_j\\ \epsilon(D)^{j}<x\end{subarray}}j^{-1}h(D)
+\frac{1}{2}(p-1)\sum_{j\geq1}\sum_{D\in\D_{p,j}^{(2)}(x)}j^{-1}h(D),\\
\hat{\pi}_{\Gamma(p)}(x^2)=&\frac{1}{2}p(p^2-1)
\sum_{j\geq1}\sum_{\begin{subarray}{c}D\in\D,p\mid u_j\\ \epsilon(D)^{j}<x\end{subarray}}j^{-1}h(D).
\end{align*}
Hence, we have
\begin{align}\label{estimate}
\sum_{j\geq1}\sum_{D\in\D_{p,j}^{(1)}(x)}j^{-1}h(D)
=&\frac{2}{p-1}\Big(\hat{\pi}_{\Gamma_1(p)}(x^2)-\frac{1}{p+1}\hat{\pi}_{\Gamma(p)}(x^2)\Big),\nt\\
\sum_{j\geq1}\sum_{D\in\D_{p,j}^{(2)}(x)}j^{-1}h(D)
=&\frac{2}{p-1}\Big(\hat{\pi}_{\Gamma_1(p)}(x^2)-\frac{1}{p}\hat{\pi}_{\Gamma(p)}(x^2)\Big).
\end{align}
Therefore, by using the prime geodesic theorem \eqref{prim} and \eqref{prim1}, 
we can obtain \eqref{d1} and \eqref{d2}.\qed

\begin{rem}
If we would get the arithmetic expressions of $Z'_{\Gamma}(s)/Z_{\Gamma}(s)$ for 
$\Gamma=\Gamma_1(p^r)\cap\Gamma(p^{r-1})$, then we could obtain more precise estimates
of $C(p)$ than ones in Corollary \ref{cor3}. 
Actually calculating the arithmetic expressions, 
we obtain 
\begin{align*}
\sum_{j\geq1}\sum_{D\in\D_{p,j}^{(2,r)}(x)}j^{-1}h(D)
=&\frac{2}{p^{3r-3}(p-1)}\Big(\hat{\pi}_{\Gamma_1(p^r)\cap\Gamma(p^{r-1}}(x^2)
-\frac{1}{p}\hat{\pi}_{\Gamma(p^r)}(x^2)\Big)\\
=&\frac{1}{p^{3r-2}}\frac{x^2}{\log{x}}+(x^{2\delta})\as x\tinf,
\end{align*}
where $\delta\in(0,1)$ is a constant and
\begin{align*}
\D_{p,j}^{(2,r)}(x):=\{D\in\D \divset  p\mmid D,p^r\nmid u_j(D),\epsilon(D)^j<x\}.
\end{align*}
Since
\begin{align*}
\bigcup_{k=1}^{r}\D_{p,j}^{(2,k)}(x)\subset \D_p(x)\subset
\bigcup_{k=1}^{r}\D_{p,j}^{(2,k)}(x)\cup\{D\in\D \divset p^r\mid u_j(D),\epsilon(D)^j<x\}
\end{align*}
and 
\begin{align*}
\sum_{j\geq1}\sum_{\begin{subarray}{c}D\in\D,\\p^r|u_j(D)\\ \epsilon(D)^j<x\end{subarray}}j^{-1}h(D)
=\frac{2}{p^{3r-2}(p^2-1)}\hat{\pi}_{\Gamma(p^r)}(x^2),
\end{align*}
we can estimate $C(p)$ as 
\begin{align*}
\sum_{k=1}^{r}\frac{1}{p^{3r-2}}\leq C(p)\leq \sum_{k=1}^{r}\frac{1}{p^{3r-2}}+\frac{1}{p^{3r-2}(p^2-1)}
\end{align*}
for any $r\geq1$.
Then $C(p)=p^2/(p^3-1)$ follows from the inequality above.
\end{rem}

\noindent 
\text{HASHIMOTO, Yasufumi}\\ 
Graduate School of Mathematics, Kyushu University.\\  6-10-1, Hakozaki, Fukuoka, 812-8581 JAPAN.\\ 
\text{hasimoto@math.kyushu-u.ac.jp}\\

\end{document}